\documentclass[12pt,reqno]{amsart}
\usepackage{setspace}
\usepackage{epsfig,hyperref}
\usepackage{amsmath, amssymb, amsthm}
\numberwithin{equation}{section}
\theoremstyle{plain}
\newtheorem{theorem}{Theorem}

\newtheorem{proposition}{Proposition}

\newtheorem{lemma}{Lemma}
\theoremstyle{remark}

\DeclareMathOperator{\supp}{supp\,}

\allowdisplaybreaks
\def\be{\begin{equation}}
\def\ee{\end{equation}}

\def\ve{\varepsilon}
\def\vp{\varphi}

\def\arrowk{^\to{\kern -6pt\topsmash k}}
\def\arrowK{^{^\to}{\kern -9pt\topsmash K}}
\def\arrowr{^\to{\kern-6pt\topsmash r}}
\def\arrowvp{^\to{\kern -8pt\topsmash\vp}}
\def\arrowf{^{^\to}{\kern -8pt f}}
\def\arrowg{^{^\to}{\kern -8pt g}}
\def\arrowu{^{^\to}a{\kern-8pt u}}
\def\arrowt{^{^\to}{\kern -6pt t}}
\def\arrowe{^{^\to}{\kern -6pt e}}
\def\tk{\tilde{\kern 1 pt\topsmash k}}
\def\barm{\bar{\kern-.2pt\bar m}}
\def\barN{\bar{\kern-1pt\bar N}}
\def\barA{\, \bar{\kern-3pt \bar A}}

\def\iint{\not \kern-4pt\int}
\begin{document}
\title{Monotone Boolean functions capture their primes}

\author{Jean Bourgain}

\address{School of Mathematics, Institute for Advanced Study, 1
Einstein Drive, Princeton, NJ 08540.}
\email{bourgain\@math.ias.edu}
\thanks{The research was partially supported by NSF grants DMS-0808042 and
DMS-0835373.}

\begin{abstract}
It is shown that monotone Boolean functions on the Boolean cube capture the expected number of primes, under the usual identification by binary 
expansion.  This answers a question posed by G.~Kalai.
\end{abstract}
\maketitle

\section
{Introduction}

Identifying the interval $\{1, 2, \ldots, 2^n-1\}$ with Boolean cube $\{0, 1\}^n$ by binary expansion, we prove the following

\begin{theorem}\label{Theorem0}
Let $f$ be a monotone Boolean function on \hfill\break
$\{0, 1\}^n \cap [x_0=1]$ and assume 
$\mathbb E[f]> c>0$, with $c$ some constant taken independent of $n$.\begin{footnote}
{A more precise quantitative version appears in (3.19)  below.}\end{footnote}
Denote $\Lambda$ the Von Mangold function. Then
\be\label{1.1}
\sum_{0< x<2^n} \Lambda(x) f(x)> \big(1-o(1)\big)2^n \mathbb E[f].
\ee
\end{theorem}

This answers an issue brought forward by G.~Kalai as part of a collection of problems related to circuit complexity, digital aspects and the Fourier-Walsh
spectrum of the Moebius and Von Mangoldt functions.  See also \cite{Gr}.

In the special case when $f$ is the majority function, we proved in \cite{B3} de-correlation from $\Lambda$, which is a stronger statement since it gives a
prime number theorem.
What is  lost in the present context of a monotone Boolean function is the invariance under the permutation group of $\{0, 1, \ldots, n-1\}$, an essential
ingredient in \cite{B3}.
In \cite{B3}, we replaced $\Lambda$ by its symmetrization $\Lambda_s$ under the permutation group, noting that $\langle f, \Lambda\rangle =\langle f,
\Lambda_s\rangle$.
This distribution $\Lambda_s$ turns out to be much better behaved than $\Lambda$, in the sense that the high order contribution of the Fourier-Walsh 
(F-W) spectrum may be handled by $L^2$-estimates.
At this point, we also invoke the tail estimate
\be\label{1.2}
\sum_{|S|>k}|\hat f(S)|^2 \lesssim k^{-\frac 12}
\ee
when $f$ is the majority function.
For general monotone Boolean functions, \eqref{1.2} has a counterpart due to Bshouty and Tamon \cite{B-T}.

\begin{proposition}
\label{Proposition1}
Let $f$ be a monotone Boolean function on $\{0, 1\}^n$.
Then for $K>1$, we have a tail estimate
\be
\label{1.3}
\sum_{|S|>K\sqrt n} |\hat f(S)|^2\lesssim K^{-1}
\ee
\end{proposition}

On the other hand, Proposition \ref{Proposition2} in the next section 
establishes \
\smallskip

\noindent
a strong bound on $\sum_{\substack{|S|<n_0\\ S\not= \phi, \{0\}}}
|\hat\Lambda(S)|^2$ for $n_0$
as large as $n^{\frac 47-\ve}$ and hence, one may hope to exploit \eqref{1.3}.
Again a direct $L^2$-approach fails and $\Lambda$ needs to be replaced by a friendlier distribution.
Keeping the monotonicity of $f$ in mind, we define
\be\label{1.4}
\tilde\Lambda =\sum_{x\in \{0,1\}^n} \Lambda(x) \sum_{j|x_j=1} \delta_{x\backslash \{j\}}
\ee
noting that
\be\label{1.5}
\langle \tilde\Lambda, f\rangle =\sum_x \Lambda(x) \sum_{j|x_j=1} f(x\backslash \{ j\}) \leq \sum_x \Lambda(x) \Big(\sum x_j\Big)f(x).
\ee
On the arithmetic side, to each prime
$$
p=2^{j_1}+2^{j_2}+\cdots+ 2^{j_k}\qquad (j_1=0)
$$
we associate the $k$ integers $x=2^{j_2} +\cdots+ 2^{j_k}, 2^{j_1} + 2^{j_3}+\cdots + 2^{j_k},\ldots, 2^{j_1}+\cdots+ 2^{j_{k-1}}$, noting that $f(x)\leq
f(p)$.

It turns out that
\be\label{1.6}
\Vert\tilde\Lambda\Vert_2 \sim\Vert\tilde\Lambda\Vert_1 \sim n
\ee
(with \eqref{1.6} expressed with normalized measure).
Further more, the  F-W coefficients $\widehat{\tilde\Lambda}(S)$ may be retrieved from the $\hat\Lambda(S)$, which
makes Proposition \ref{Proposition2} also applicable to $\tilde\Lambda$.
\medskip

In \S2, we prove Proposition \ref {Proposition2}, a result of independent interest.
The argument is in fact quite similar in many aspects to the analysis 
in \cite{B1} and we only indicate the main points.
\medskip

In \S3, the function $\tilde\Lambda$ is analyzed and the scheme described above worked out; inequality \eqref{1.1} is established in a more precise form.

As pointed out earlier, \eqref{1.1} does not provide an asymptotic formula and it may be unreasonable to expect one in this generality.

\medskip
\noindent
{\bf Acknowledgement:} The author is grateful for various communications with G.~Kalai on the subject.

\section
{On the Fourier-Walsh coefficients of the von Mangoldt function}

Let $N=2^n$ and identify $\{0, 1, \ldots, N-1\}$ with $\{0, 1\}^n$ by binary expansion
\be\label{2.1}
x=\sum_{0\leq j< n} x_j2^j \ \text { with } \ x_j=0, 1.
\ee
The Walsh system $\big \{w_S; S\subset \{0, 1, \ldots, n-1\}\big\}$ is defined by
\be\label{2.2}
w_S(x)=\prod_{j\in S}\ve_j\  \text { with } \ \ve_j= 1-2x_j \in \{1, -1\}.
\ee
If $f$ is a function on $\{0, 1\}^n$, we have
\be\label{2.3}
f(x) =\sum_S \hat f(S) w_S(x) \ \text { where } \ \hat f(S) =2^{-n}\sum_{x\in\{0, 1\}^n} f(x) w_S(x).
\ee

Considering the restriction of the Von Mangoldt function $\Lambda$ to \hfill\break
$\{1, \ldots, N-1\}$ as a function on the Boolean cube $\{0, 1\}^n$, the following
estimate holds on its Fourier-Walsh spectrum $\{\hat\Lambda (S)\}$.

\begin{proposition}\label{Proposition2}
Let $n_0< n^{4/7} (\log n)^{-2}$.
Then for $n$ large enough
\be\label{2.4}
\sum_{|S|\leq n_0, S\not= \phi, \{0\}} |\hat\Lambda(S)|^2 < e^{-n^{3/7}}.
\ee
\end{proposition}

The condition $S\not= \{0\}$ of course is parity related.

Note that for application in conjunction with Proposition 1, it is essential that $n_0$ brakes the $n^{1/2}$-barrier (cf. the related work by
B.~Green \cite{Gr} about correlation with $AC(0)$ circuits).

A similar statement was established in \cite{B1} for the Moebius function $\mu$ (stated with a worse estimate).
Despite many similarities, there are significant differences in the analysis and we rather follow the treatment in \cite{B2}.
The problem studied in \cite{B2} is the evaluation of
\be\label{2.5}
\sum^N_1 \Lambda(x) f(x)
\ee
with $f$ of the form $f= 1_{[x_{j_1}=\alpha_1, \ldots, x_{j_r}=\alpha_r]}$ (i.e. primes with certain prescribed binary digits).
Most of the analysis in \cite {B2} is presented in the context of a general function $f$ however.
The following technical statement summarizes \S 1, 2, 3 from \cite{B2}.

\begin{lemma}\label{Lemma1}
Assume $f$ a function on $\{1, \ldots, N-1\}$ supported on the odd integers and of mean zero.  Then
$$
\Big|\sum^N_1 \Lambda (x) f(x) \Big|\lesssim
$$
\be\label{2.6}
\frac N{\sqrt B} n^3\Vert\hat f\Vert_1
\ee
\be\label{2.7}
+ B\Vert\hat f\Vert_1 \Big[\max_{\substack {u<N\\ \mathcal X\in\mathcal G}} |\psi(u, \mathcal X)|\Big]
\ee
\be\label{2.8}
+n^2\Big[\sum_{q<B, q\, sf} \ \frac {\kappa (q)} q\Big] \Vert f\Vert_1
\ee
\be\label {2.9}
+ B\exp \big(-c n^{3/5} (\log n)^{-\frac 15}\big) \Vert f\Vert_1
\ee
\be\label{2.10}
+ n^5\Big\{\sum_{\mathcal X_1\in \mathcal B}\nolimits ^* \frac {q_1}{\phi (q_1)} \Big[\sum_{q_3< B, q_3 s f,  \text { odd}} \ \frac{\alpha (q_1, q_3)}{q_3}\Big]\Big\}
\Vert f\Vert_1.
\ee

Here $\hat f$ refers to the usual Fourier transform, $\psi(u, \mathcal X) =\sum_{x<u} \Lambda(x) \mathcal X(x)$ and $B$ is a parameter.
The classes $\mathcal G$ and $\mathcal B$ corresponding to a subdivision of the Dirichlet characters according to the zeros of their corresponding $L$-function.
More specifically, let $T>B, \log T\sim \log B$ be another parameter and denote
\be\label{2.11}
\eta(\mathcal X) =\min \{ 1-\beta; \rho=\beta+i\gamma, |\gamma|<T \text { zero of } L(s, \mathcal X)\}.
\ee
Then
\be\label{2.12}
\mathcal G=\Big\{\mathcal X (\text{mod\,} q) \text { non-principal}; q<B \text{ and } \eta(\mathcal X)> C\frac {\log T}n\Big\}
\ee
and
\be\label{2.13}
\mathcal B=\big\{ \mathcal X(\text{mod\,} q) \text { non-principal} ; q< B\big\}\backslash \mathcal G.
\ee
In \eqref{2.10}, $\sum^*$ refers to summation over primitive characters.

Remains to specify $\kappa(q), \alpha (q_1, q_3)$.
They are required to satisfy the respective inequalities
\be\label{2.14}
\Big|\sum_{x<N, q|x} f(x)\Big| \leq \frac{\kappa(q)}q \Vert f\Vert_1 \text { for $q<B$ odd and square-free}
\ee
and
\be\label{2.15}
\Big|\sum_{x\in J, q_3|x} f(x) \mathcal X_1(x)\Big| \leq \frac {\alpha(q_1, q_3)}{q_3} 
\ \sum_{x\in J} |f(x)|
\ee 
for $q_1, q_3 <B, (q_1, q_3)=1$, $q_3$ odd and square-free, $\mathcal X_1(\text{mod\,} q_1)$ primitive and $J\subset [1, N]$ an arbitrary interval of size $\sim
\frac NB$.
\end{lemma}

The above the statement is established using the circle method.  A few words of explanation about the different contributions.
The term \eqref{2.6} is the minor arcs contribution, estimated using Vinogradov's inequality.
The major arcs are analyzed the usual way, using Dirichlet characters; \eqref{2.8}, \eqref{2.9} account for the contribution of the principal characters.

The contribution of the non-principal characters $\mathcal X\in \mathcal G$ is expressed by \eqref{2.7} and those in $\mathcal B$ by \eqref{2.10}.
It is important to note that the savings in \eqref{2.8}, \eqref{2.10} depend on $\kappa(q)$ and $\alpha (q_1, q_3)$ in inequalities \eqref{2.14}, \eqref{2.15}.

Following \cite{B2}, \S6, take
\be\label{2.16}
\log T\sim n^{4/7} (\log n)^{-3/7}
\ee
which in \eqref {2.7} gives an estimate
\be\label{2.17}
|\psi(u, \mathcal X)|\lesssim \frac {n^2N}{T} \ \text { for } \mathcal X\in\mathcal G, u<N.
\ee
In order to deduce Proposition \ref{Proposition1} from Lemma \ref{Lemma1}, the natural choice for $f$ would be the truncation
\be\label{2.18}
f=\Big[\sum_{\substack {S\subset\{1, \ldots, n-1\}\\ 0<|S|<n_0}} \hat\Lambda (S) w_S\Big] (1-\ve_0) =\sum \hat f(S) w_S.
\ee
Similarly to \cite{B2}, we make the following modification of $f$.

Let $m\in \mathbb Z_+$ be another parameter, satisfying
\be\label{2.19}
\log B<m, n\sim\log B
\ee
and partition $\{0, 1, \ldots, n-1\}$ in intervals $J_\alpha$ of size $\frac nm$.

Let
\be\label{2.20}
K_0= 1+ 10\frac {n_0m}n
\ee
and define
\be\label{2.21}
\omega_\alpha(S) =\begin{cases} 1\text { if } |S\cap J_\alpha|\geq K_0\\ 0 \text { otherwise.}\end{cases}
\ee
Hence $\sum_\alpha \omega_\alpha (S) \leq\frac {|S|}{K_0}$ and by \eqref{2.20}
$$
\sum_\alpha \Big[\sum_{|S|<n_0}|\hat\Lambda(S)|^2 \big(1-\omega_\alpha (S)\big)\Big] \geq \frac n{2m} \sum_{|S|<n_0} |\hat\Lambda(S)|^2.
$$
We may therefore redefine for some $\alpha =1, \ldots, \big[\frac nm\big]$ the function $f$ as
\be\label{2.22}
f=f_\alpha =\sum_{\substack {|S|\leq n_0\\ |S\cap J_\alpha|<K_0}} \hat f(S) w_S
\ee 
and still satisfy
\be\label{2.23}
\sum\hat\Lambda(S) \hat f(S) \gtrsim \sum_{|S|<n_0} |\hat\Lambda (S)|^2.
\ee
Note that $\Vert f\Vert_2 \leq \Vert\Lambda\Vert_2 \leq \sqrt{Nn}, \Vert f\Vert_1 \leq \sqrt nN$ and, cf. \cite{B1}, \cite{B2}
\be\label{2.24}
\Vert\hat f\Vert_1< (Cn)^{n_0}.
\ee
Recalling \eqref{2.16}, \eqref{2.17}, it follows that
\be\label{2.25}
\eqref{2.6}+\eqref{2.7}+\eqref{2.9} < (Cn)^{n_0} N\Big(\frac 1{\sqrt B}+\frac BT\Big)<(Cn)^{n_0}\frac N{\sqrt B}
\ee
provided $B< \sqrt T$.

Obtaining the required bounds on
\be\label{2.26}
\sum_{\substack{q< B\\ q \, sf \text { and odd}}} \ \frac{\kappa (q)}{q}
\ee
and
\be\label{2.27}
\sum_{\substack{q_3<B\\ q_3 \, sf \text{ and odd}}} \ \frac {\alpha(q_1, q_3)}{q_3}
\ee
in \eqref{2.8}, \eqref{2.10} is a non-trivial and crucial input in \cite{B1}, \cite{B2}.
Essential use is made in the argument of the property $|S\cap J_\alpha|<K_0$ in \eqref{2.22}.
The basic idea is as follows.
Assume that  $|S\cap J_\alpha|< K_0$ if $\hat f(S)\not= 0$.
Then the function $f$ restricted to a translate of a progression of the form $2^j.I$, $I=\{1, \ldots, m\}$, has a low order Walsh expansion.

Switching to an expansion in the trigonometric system (with suitable approximation of the Walsh functions) permits then to analyze the l.h.s. in
\eqref {2.14}, \eqref{2.15}.

This analysis was carried out in \cite{B1}, \S4, 5 for a function $f$ of the form $1_{[x_{j_1}=\alpha_1, \ldots, x_{j_r}=\alpha_r]}$.
Some adjustments of this rather tedious argument (noting also that $\mathbb E[f]=0$ here, by assumption) permit us to obtain an estimate of the form
\be\label{2.28}
\eqref{2.26}, \eqref{2.27} < e^{-c\frac n{n_0}}.
\ee
Note that in \cite{B2}, \S4, this estimate on \eqref{2.27} is proven for a function $f$ of the type considered above.
\medskip

From \eqref{2.25}, \eqref{2.29}, it follows that
\be\label{2.29}
\Big|\sum^N_{1} \Lambda(x) f(x)\Big| < N \big((Cn)^{n_0} B^{-\frac 12} + e^{-c\frac n{n_0}} (n^3+|\mathcal B|)\big)
\ee
where the size of $\mathcal B$ is controlled from classical zero-density estimates
\be\label{2.30}
|\mathcal B|\lesssim (TB)^{8\eta_*}, \eta_*= C\frac {\log T}n
\ee
(cf. \cite{B1}, \S6).
Taking, according to \eqref{2.16}
\be\label{2.31}
\log B=c\min \Big(\frac n{\sqrt {n_0}}, \frac {n^{4/7}}{(\log n)^{3/7}}\Big)
\ee
and recalling the assumption on $n_0$ in Proposition \ref{Proposition1}, gives
\be\label {2.32}
\eqref{2.29} < \exp \Big\{ -c\min\Big[\frac n {n_0}, \frac {n^{4/7}}{(\log n)^{3/7}}\Big]\Big\}.
\ee
Inequality \eqref{2.4} follows.

\section
{The function $\tilde\Lambda$ and proof of the theorem}

Recall that $\tilde\Lambda$ is defined as
\be\label{3.1}
\tilde\Lambda= \sum_x \Lambda(x) \sum_{j|x_j=1} \delta_{x\backslash \{j\}}
\ee
with $\delta_{x\backslash \{j\}}$ the Dirac at $x\backslash \{j\}$.
From \eqref{1.5}, implied by monotonicity,
\be\label{3.2}
\langle \tilde\Lambda, f\rangle \leq n\langle \Lambda, f\rangle.
\ee
We start by evaluating $\Vert\tilde\Lambda\Vert_2$.
Clearly, by \eqref{3.1}
\be
\begin{aligned}\label{3.3}
\sum_{y\leq N}\tilde \Lambda(y)^2 &\leq \sum_{y\leq N} \ \sum_{0\leq j, k<n} \Lambda(y+2^j)\Lambda (y+2^k)\\
&=\sum_{0\leq j, k<n} \ \sum_{x<2N}\Lambda(x) \Lambda(x+2^k-2^j).
\end{aligned}
\ee
Using a simple upper bound sieve (cf. \cite{I-K}, Th. 6.7), one has for $j\not= k$
\be\label {3.4}
\sum_{x<2N} \Lambda(x) \Lambda (x+2^k-2^j)\lesssim N
\ee
implying that
\be\label{3.5}
\Vert\tilde\Lambda
\Vert_2\lesssim n\sqrt N.
\ee
\medskip

Next, we compute the F-W coefficient of $\tilde \Lambda$.

\medskip
Rewrite \eqref{3.1} as
$$
\tilde\Lambda =\sum_x\Lambda (x) \Big[ \delta_{(0, x_1, \ldots, x_{n-1})} +\delta_{(x_0, 0, x_2, \ldots, x_{n-1})} +\cdots+\delta_{(x_0, \ldots, x_{n-2},0)}
-\Big (n-\sum x_j\Big)\delta_x\Big].
$$

For $S\subset\{0, 1, \ldots, n-1\}$, we get
\be\label{3.6}
\sum_x \Lambda(x) w_S(0, x_1, \ldots, x_{n-1})=\begin{cases} N \hat\Lambda (S)  \ \text { if } 0\not\in S\\ N\hat \Lambda(S\backslash \{0\}) \ 
\text { if } 0\in S.\end{cases}
\ee
Also from the formulas $\ve_j=1- 2 x_j$, $n-\sum\limits_0^{n-1} x_j =\sum\limits_j^{n-1} \frac{1+\ve_j}2$, we obtain
\be\label{3.7}
\begin{aligned}
&\sum_x\Lambda(x) (n-\sum x_j)w_S(x) = N\langle\Lambda, \Big(\sum_0^{n-1} \frac {1+\ve_j} 2 \Big) \prod_{j\in S} \ve_j\rangle=\\
&N\Big\{\frac n2\hat\Lambda (S)+\frac 12 \sum_{j\in S} \hat\Lambda (S\backslash \{j\})+\frac 12\sum_{j\not\in S}\hat\Lambda(S\cup\{j\})\Big\}.
\end{aligned}
\ee
Hence
\be\label{3.8}
\widehat{\tilde\Lambda}(S) = \Big(\frac n 2-|S|\Big) = \tilde\Lambda (S)+\frac 12 \sum_{j\in S} \hat\Lambda (S\backslash \{j\}) -\frac 12\sum_{j\not \in S}
\hat\Lambda (S\cup \{j\}).
\ee

In particular, it follows that
\be\label{3.9}
\mathbb E[\tilde\Lambda]= \frac n2 -\frac 12\sum_j\hat\Lambda(\{j\}) = \frac {n-1}2+o(1)
\ee
\be\label{3.10}
\Vert\tilde\Lambda\Vert_1 \sim nN.
\ee
Combined with \eqref{3.5}, it follows that the distribution $\tilde\Lambda$ is essentially `flat'.

Write using F-W expansion and \eqref{3.9}
\be\label{3.11}
\langle f, \tilde\Lambda\rangle =\mathbb E[f] \,  \mathbb E[\tilde\Lambda]+ (3.12) + (3.13)
> \Big(\frac n2-1\Big) \mathbb E[f] + (3.12)+ (3.13) =\eqref{3.11}
\ee
with
$$
(3.12) =\sum_{\substack {S\not= \phi\\ |S|<K\sqrt n}}\hat f(S) \hat{\tilde\Lambda} (S)
$$
and
$$
(3.13) =\sum_{|S|\geq K\sqrt n} \hat f(S) \hat{\tilde \Lambda} (S)
$$
($K$ a large parameter).

In view of \eqref{1.3}, \eqref{3.5}
$$
(3.13) \lesssim K^{-\frac 12} \frac {\Vert\tilde\Lambda \Vert_2}{\sqrt N} \lesssim K^{-\frac 12} n.
\eqno{(3.14)}
$$
It follows from \eqref{3.8} and Proposition \ref{Proposition2} that certainly
$$
\sum_{2<|S|\leq K\sqrt n} |\hat{\tilde \Lambda} (S)|^2 < n^2 e^{-n^{3/7}}
\eqno{(3.15)}
$$
assuming
$$
K<n^{\frac 1{14}} (\log n)^{-2}.
\eqno{(3.16)}
$$
For $S=\{0\}$, by \eqref{3.8}, (3.15)
$$
\begin{aligned}
\hat{\tilde \Lambda}(\{0\}) &=\Big(\frac n2-1\Big) \hat\Lambda(\{0\})+\frac 12\hat\Lambda(\phi)+O(e^{-\frac 13 n^{3/7}})\\
&= \frac{3-n}2+ O\big(e^{-\frac 13 ^{n^{3/7}}} \big)
\end{aligned}
$$
since $\hat\Lambda(\{0\}) =\mathbb E[\Lambda\ve_0]= -\hat\Lambda(\phi)$ by parity.
\medskip

If $S=\{j\}, 0<j<n$,
$$
\hat{\tilde\Lambda} (\{j\})=\frac 12\hat{\tilde\Lambda}(\phi)+ O(e^{-\frac 13 n^{3/7}})=\frac 12+ O(e^{-\frac 13 n^{3/7}} ).
$$
For $S=\{0, j\}$,
$$
\hat{\tilde\Lambda} (S)=\frac 12 \hat\Lambda(\{0\})+ O(e^{-\frac 13 n^{3/7}}) =-\frac 12+O(e^{-\frac 13 n^{3/7}})
$$
and $|\hat{\tilde\Lambda}(S)|< O(e^{-\frac 13 n^{3/7}})$ if $|S|=2, 0\not\in S$.
\medskip

Therefore
$$
\begin{aligned}
(3.12) &= -\frac {n-3} 2\hat f(\{0\})+\frac 12 \sum_{0<j<n}\hat f(\{j\})-\frac 12 \sum_{0<j< n} \hat f(\{0, j\})+O( e^{-\frac 13 n^{3/7}})\\
&=-\frac {n-3} 2 \mathbb E[f.\ve_0]+ \frac 12\sum_{0<j<n}\mathbb E[(1-\ve_0)\ve_j f]+ O(e^{-\frac 13 n^{3/7}})\\
&=\frac {n-3}2 \mathbb E[f]+\sum_{0<j<n} \mathbb E[\ve_j f]+O(e^{-\frac 13 n^{3/7}})\\
&=\frac {n-3}2\mathbb E[f] +O(\sqrt n)
\end{aligned}
\eqno{(3.17)}
$$
again using that $\supp f\subset[x_0=1]$.
\medskip

Substituting (3.14), (3.17) into \eqref{3.11},
$$
\langle f, \tilde\Lambda\rangle > \big(n-O(1)\big)\mathbb E[f]+O(\sqrt n)+O(K^{-\frac 12} n).
\eqno{(3.18)}
$$
From (3.16), \eqref{3.2}
$$
\langle f, \Lambda\rangle \geq \mathbb E[f]+O(n^{-\frac 1{28}}\log n).
\eqno{(3.19)}
$$
This gives \eqref{1.1}, with in fact $n^{-\frac 1{28}} \log n$ in place of $o(1)$.


\begin{thebibliography}
{xxxxxxxx}

\bibitem[B1]{B1} J.~Bourgain, \emph{Prescribing the binary digits of the primes}, to appear in Israel J.~Math.

\bibitem[B2]{B2} J.~Bourgain, \emph{On the Fourier-Walsh spectrum of the Moebius function}, to appear in Israel J. Math.

\bibitem[B3]{B3} J.~Bourgain, \emph {A prime number theorem for the majority function}, preprint.

\bibitem[B-T]{B-T} N.~Bshouty, C.~Tamon, \emph{On the Fourier spectrum of monotone functions}, Journal of the ACM, Vol. 43, N4 (1996), 747--770.

\bibitem[Gr]{Gr} B.~Green, \emph{On (not) computing the Moebius function using bounded depth circuits}, ArXiv: 1103,4991.

\bibitem[Ka]{Ka} G.~Kalai, \emph{Private communications}

\bibitem[I-K]{I-K} H.~Iwaniec, E.~Kowalski, Analytic Number theory, AMS.
\end{thebibliography}
\end{document}